\documentclass[11pt,reqno]{amsart}
\usepackage{color}
\usepackage[active]{srcltx}
\usepackage{a4wide}
\usepackage{amssymb, amsmath, mathtools}
\usepackage{graphicx}
\usepackage{float}
\usepackage[active]{srcltx}
\usepackage[ansinew]{inputenc}
\usepackage{hyperref}

\theoremstyle{plain}
\newtheorem{theorem}{Theorem}[section]

\newtheorem{lemma}[theorem]{Lemma}

\newtheorem{corollary}[theorem]{Corollary}

\theoremstyle{remark}

\newtheorem{remark}[theorem]{Remark}

\numberwithin{equation}{section}

\newcommand{\NN}{{\mathbb{N}}}
\newcommand{\ZZ}{{\mathbb{Z}}}
\newcommand{\RR}{{\mathbb{R}}}
\newcommand{\EU}{{\mathbb{S}}}

\newcommand{\Fix}{{\text{Fix}}}

\newcommand{\dpt}{\displaystyle}

\begin{document}

\title[Torus-breakdown near a Bykov attractor]{Torus-breakdown near a Bykov attractor: \\  a case study}

\author[Lu\'isa Castro]{Lu\'isa Castro \\ Center for Health Technology and Services Research - CINTESIS, \\ University of Porto, \\
Rua Dr. Pl\'acido da Costa, 4200-450 Porto, Portugal \\ \\ } 
\address[L. Castro]{Center for Health Technology and Services Research - CINTESIS, University of Porto, Rua Dr. Pl\'acido da Costa, 4200-450 Porto, Portugal}
\email[L. Castro]{luisacastro@med.up.pt} 

\bigbreak

\author[Alexandre A. P. Rodrigues]{\\ \\ Alexandre A. P. Rodrigues \\ Centro de Matem\'atica da Universidade do Porto \\ 
and Faculdade de Ci\^encias da Universidade do Porto\\
Rua do Campo Alegre 687, 4169--007 Porto, Portugal\\}
\address[A. A. P. Rodrigues]{Centro de Matem\'atica da Universidade do Porto\\
and Faculdade de Ci\^encias da Universidade do Porto\\ \\
Rua do Campo Alegre 687, 4169--007 Porto, Portugal}
\email[A.A.P.Rodrigues]{alexandre.rodrigues@fc.up.pt}


\begin{abstract} 
There are few explicit examples in the literature of vector fields exhibiting complex
dynamics that may be proved analytically. 
This paper reports numerical experiments performed for an explicit two-parameter family of vector fields unfolding an attracting heteroclinic network, linking two saddle-foci with $ (\mathbb{SO}(2) \oplus \ZZ_2)$--symmetry. 
The vector field is the restriction to $\EU^3$ of a polynomial vector field in $\RR^4$.
We investigate global bifurcations due to symmetry-breaking  and we detect strange attractors via a phenomenon called \emph{Torus-Breakdown theory}. We explain how an attracting torus gets destroyed by following the changes in the invariant manifolds of the saddle-foci. 

  Although a complete understanding of the corresponding bifurcation diagram and the mechanisms underlying the dynamical changes is still out of reach, using a combination of theoretical tools and computer simulations, we have uncovered some complex patterns for the symmetric family under analysis. This also suggests a route to obtain rotational horseshoes; additionally, we give an attempt to elucidate some of the bifurcations involved in an Arnold wedge. 
\end{abstract}

\maketitle
\setcounter{tocdepth}{1}
\textbf{
We explicitly construct  a two-parameter family of polynomial differential equations, in which each parameter controls a type of symmetry-breaking. 
We discuss global bifurcations that occur as the two parameters vary, namely the emergence of strange attractors.  This  route  to chaos has been a recurrent concern on nonlinear dynamics during the last decades. }

\textbf{
We use a systematic method to construct examples of vector fields with simple forms that make their dynamical properties amenable to analytic proof. The method consists of using symmetry to obtain ``gradient-like dynamics'' and then choosing special symmetry-breaking nonlinear terms with a simple form that preserve the required properties and introduce some desired behaviour.}

\textbf{
For this class of examples, we show  the existence of many complicated dynamical objects and the complex transition between different types of dynamics, ranging from an attracting torus to H\'enon-like strange attractors, as a consequence of the \emph{Torus-Breakdown theory}.  Different symmetries are broken gradually and independently; different phenomena are associated to the inclusion of different symmetry-breaking terms.}

\section{Introduction}\label{intro}

A compact attractor is said \emph{strange} if it contains a dense orbit with at least one positive Lyapunov exponent. 
 The rigorous proof of the strange character of an invariant set is a great challenge and the proof of their prevalence in the space of parameters (persistence with respect to the Lebesgue measure) is a very involving task. 

In the present paper, illustrating theoretical results of \cite{Rodrigues2019}, we construct a two-parameter symmetric family of polynomial differential equations inspired in \cite {Aguiar_tese}, and we show some evidences for  the  existence of strange attractors in its unfolding. Using the \emph{maximal Lyapunov exponent} along an orbit, we explore numerically a mechanism to  
obtain strange attractors in the unfolding of an equivariant vector field. The abundance of strange attractors will be a consequence of the \emph{Torus-Breakdown Theory} developed in \cite{AS91, AHL2001, AH2002, Aronson}.

\subsection{Lyapunov exponents}

 A \emph{Lyapunov exponent} associated to a solution of a differential equation is an average exponential rate of divergence or convergence of nearby trajectories in the phase space. As they measure the rate at which the dynamics creates or destroys information, the Lyapunov exponents equal in number the dimension of the phase space and allow us to distinguish between chaos and regular dynamics (stable periodicity).

 A positive exponent reflects the existence of a direction in which the system experiences the repeated stretching and folding that mixes nearby states on the attractor. Thus, the long-term behavior of an initial condition with a positive Lyapunov exponent cannot be predicted, and this is one of the most common features of chaos. This is the key idea behind several numerical experiments with chaotic dynamical systems. Since nearby solutions may correspond to numerical almost identical states, the presence of exponential orbital divergence implies that trajectories whose initial conditions are hard to distinguish will soon depart, and most likely behave afterwards quite differently. 
 
 The study of the number of positive Lyapunov exponents along well chosen orbits\footnote{In general, these well chosen orbits correspond to the unstable manifolds of invariant saddles.} motivates one of the most powerful computational techniques available to build a reliable approximation of a bifurcation diagram.

\subsection{This article}
We start  with a detailed study of the  two-parameter polynomial differential equation: equilibria, symmetries, flow-invariant sets, relative positions of the invariant manifolds, heteroclinic connections, Lyapunov stability. Then, 
we perform several illustrative computer experiments using  Matlab (R219b, Mathworks, Natick, MA, USA)  for this family of vector fields. The periodic or chaotic nature of solutions could only be determined on a case-by-case examination by fixing parameters and investigating the dynamics for well chosen initial conditions. Additional care was needed while interpreting the numerical integration of these flows since, for some parameters, they exhibit \emph{quasi-stochastic attractors} \cite{Afraimovich83} and these are prone to rounding errors that may ruin the simulations. 
\medbreak

This article is organised as follows. In Section \ref{s:setting}, based on \cite{Aguiar_tese, Rodrigues2019} we revisit the setting of symmetric Bykov attractors and the main theoretical results that state the existence of strange attractors. In Section \ref{s:example1}, we exhibit an explicit two-parameter family of equivariant vector fields that will be the object of consideration throughout the paper. The construction of the vector field is amenable to the analytic proof of the features that guarantee complex behaviour. We illustrate   dynamical phenomena in that example (when parameters vary) going from an attracting torus to horseshoes. In Section \ref{interpretation}, we describe the expected theory about the emergence of strange attractors from an attracting torus breaking. We will see that the numerics of this section agree perfectly well with the existing theory on the topic. Finally, in Section \ref{s:discussion}, we discuss the results relating them with others in the literature. For reader's convenience, we have compiled at the end of the manuscript a list of definitions in a short glossary.

\section{The theory -- an overview}
\label{s:setting}
Our object of study is the dynamics around an attracting  heteroclinic network for which we give a rigorous description here. 
For each subset $M\subset \EU^3$, we denote by $\overline{M}$ its topological closure in $\EU^3$. In order not to interrupt the flow of ideas, we refer to Appendix \ref{Definitions} for the technical definitions of some of the terms.
\subsection{The organising center}
For $\varepsilon>0$ small enough, consider the two-parameter family of $C^3$-smooth differential equations
\begin{equation}
\label{general2}
\dot{x}=f_{(A, \lambda)}(x)\qquad x\in \EU^3 \qquad A, \lambda \in [0, \varepsilon] 
\end{equation}
and denote by $\varphi_{(A, \lambda)}(t,x)$, $t \in \RR$, the associated flow, satisfying the following hypotheses for $A=0$ and $\lambda=0$:

\bigbreak
\begin{enumerate}
 \item[\textbf{(P1)}] \label{B1}  There are two hyperbolic equilibria, say $O_1$ and $O_2$.
 \bigbreak
 \item[\textbf{(P2)}] \label{B2} The spectrum of $Df_X$ is:
 \medbreak
 \begin{enumerate}
 \item[\textbf{(P2a)}] $E_1$ and $ -C_1\pm \omega_1 i $ where $C_1>E_1, \, \, \omega_1>0$, \quad  for $X=O_1$;
 \medbreak
 \item[\textbf{(P2b)}] $-C_2$ and $ E_2\pm \omega_2 i $ where $C_2> E_2,\, \,  \omega_2>0$,  \quad for $X=O_2$.
 \end{enumerate}
 \end{enumerate}
\bigbreak 
 Thus the equilibrium $O_1$ possesses a 2-dimen\-sional stable and $1$-dimen\-sional unstable manifold and the equilibrium $O_2$ possesses a 1-dimen\-sional stable and $2$-dimen\-sional unstable manifold. We assume that:
 \begin{enumerate}
 \bigbreak
  \item[\textbf{(P3)}]\label{B3} The sets $\overline{W^u(O_2)}$ and $\overline{W^s(O_1)}$ coincide and  $\overline{W^u(O_2) \cap W^s(O_1)}$ consists of a two-sphere (also called the $2D$-connection) containing $O_1$ and $O_2$. 
  \end{enumerate} \bigbreak
  and 
   \bigbreak
   \begin{enumerate}
\item[\textbf{(P4)}]\label{B4} There are two trajectories, say  $\gamma_1, \gamma_2$, contained in  $W^u(O_1)\cap W^s(O_2)$,  one in each connected component of $\EU^3\backslash \overline{W^u(O_2)}$ (called the $1D$-connections).
\end{enumerate}

\bigbreak
The two equilibria $O_1$ and $O_2$, the two-dimensional heteroclinic connection from $O_2$ to $O_1$ refered in \textbf{(P3)} and the two trajectories listed in  \textbf{(P4)}  build a heteroclinic network we will denote hereafter by $\Gamma$. 
This set has a \emph{global attracting} character, this is why it will be called by \emph{Bykov attractor}; terminology and details in \eqref{attracting set} and \eqref{ss:Bylov_cycle}.  In particular, we may find an open neighborhood $\mathcal{U}$ of the heteroclinic network $\Gamma$ having its boundary transverse to the flow associated to the vector field $f_{(0,0)}$ and such that every solution starting in $\mathcal{U}$ remains in it for all positive time and is forward asymptotic to $\Gamma$ (Lemma 2.1 of \cite{Rodrigues2019}). 
\medbreak
There are two possibilities for the geometry of the flow around each saddle-focus of the network $\Gamma$, depending on the direction the solutions turn around $[O_1 \rightarrow O_2]$. We assume that:
\medbreak
\begin{enumerate}
\item[\textbf{(P5)}] \label{B5} The saddle-foci $O_1$ and $O_2$ have the same chirality (details in \eqref{chirality_def}).
\end{enumerate}

\bigbreak

For  $r \geq 3$, denote by  $\mathfrak{X}^r(\EU^3)$, the set of two-parameter families of $C^r$--vector fields on $\EU^3$ endowed with the $C^r$--Whitney topology, satisfying Properties \textbf{(P1)--(P5)}.

\subsection{Perturbing terms}
With respect to the effect of the two parameters $A$ and $\lambda$ on the dynamics, we assume that:
 \medbreak
\begin{enumerate}
\item[\textbf{(P6)}] \label{B5}  For all $A> \lambda\geq 0$, the two trajectories within $W^u(O_1)\cap W^s(O_2)$ persist.
\end{enumerate}
 \medbreak
\begin{enumerate}
\item[\textbf{(P7)}]\label{B7} For all $A> \lambda\geq 0$, the two-dimensional manifolds $W^u(O_2)$ and $W^s(O_1)$ do not intersect.
\end{enumerate}

 \medbreak

Rodrigues \cite{Rodrigues2019} created a model assuming an extra technical hypothesis (for $A>~\lambda \geq 0$):
\medbreak

\begin{enumerate}
\item[\textbf{(P8)}] \label{B8} The transitions along the connections $[O_1 \rightarrow  O_2]$ and $[O_2 \rightarrow  O_1]$ are given, in local coordinates, by the Identity map and, up to high order terms, by $$(x,y)\mapsto (x,y+A + \lambda \Phi(x))$$ respectively, where $\Phi:\EU^1 \rightarrow \EU^1$ is a Morse smooth function with at least two non-degenerate critical points ($\EU^1=\RR \pmod{2\pi}$).
\end{enumerate}
\medbreak


\medbreak
Hypothesis \textbf{(P8)} is natural when we consider the \emph{Melnikov integral} \cite{GH} applied to a  differential equation of the type \eqref{general2}. The distance between $W^u_{\text{loc}}(O_2)$ and $W^s_{\text{loc}}(O_1)$ in a given cross section to $\Gamma$ may depend on a variable $x$ and it decomposes as 
$$
Mel(x)= Mel_1 + \lambda \, Mel_2 (x)
$$
where $Mel_1 \equiv A $ gives the averaged distance between $W^u_{\text{loc}} (O_2)$ and $W^s_{\text{loc}} (O_1)$  in the given cross section    and $Mel_2\equiv \Phi$ describes fluctuations of the unstable manifold of $O_2$. To simplify the notation, in what follows we will sometimes drop the subscript $(A, \lambda)$, unless there is some risk of misunderstanding.

\subsection{Notation}
\label{notation1}
From now on, we settle the following notation:
\begin{equation}
\label{constants}
\delta_1 = \frac{C_1}{E_1 }>1, \qquad \delta_2 = \frac{C_2}{E_2 }>1, \qquad \delta=\delta_1\, \delta_2>1 
 \end{equation}
and
\begin{equation}
\label{constants2}
   K_\omega= \frac{E_2 \, \omega_1+C_1\, \omega_2 }{E_1E_2}>0.
 \end{equation}
The constant $K_\omega$ will be called the \emph{twisting number} of $\Gamma$.  From now on, denote by $\mathfrak{X}^r_{\text{Byk}}(\EU^3)$ the set of two-parameter families of $C^r$--vector fields that satisfy the conditions  \textbf{(P1)}--\textbf{(P8)}.  The parameters $A$ and $\lambda$ are supposed to be small.

\subsection{The results}

According to \cite{Rodrigues2019}, we may draw, in the first quadrant,  two smooth curves, the graphs of $h_1$ and $h_2$, such that:
\begin{enumerate}
\medbreak
\item  $\dpt h_1(K_\omega)=\frac{1}{\sqrt{1+K_\omega^2}}$ and $ \dpt h_2(K_\omega)=  \frac{\exp\left(\frac{6\pi}{K_\omega \, }\right)-1}{\exp\left(\frac{6\pi}{K_\omega \, }\right)-1/6}$;
\medbreak
\item the region below the graph of $h_1$  corresponds to flows having an invariant and attracting torus with \emph{zero topological entropy}  (regular dynamics);
\medbreak
\item the region above the graph of $h_2$ corresponds to vector fields whose flows exhibit \emph{rotational horseshoes} in the sense of Passegi \emph{et al} \cite{Passeggi} -- see \eqref{ss: rotational horseshoe}.
\end{enumerate}
\medbreak
Under some conditions on the parameters and on the eigenvalues of the linearisation of the vector field at the saddle-foci, the author of  \cite{Rodrigues2019} proved the existence of  H\'enon-like strange attractors near the ``\emph{ghost}'' of the Bykov attractor:

\medbreak

\begin{theorem}[\cite{Rodrigues2019}, adapted]
\label{Rod_th}
Let $f_{(A, \lambda)} \in\mathfrak{X}_{\emph{Byk}}^3(\EU^3)$. 
 Fix $K_\omega^0>0$.  In the bifurcation diagram $\left(A,\frac{\lambda}{A} \right)$, where $(A, \lambda)$ is  such that $ h_1(K_\omega^0) <\frac{\lambda}{A} <h_2(K_\omega^0)$, there exists a positive measure set  $\Delta$ of parameter values, so that for every $\lambda/A\in \Delta$, the flow of \eqref{general2}  admits a strange attractor  of H\'enon-type with an ergodic SRB measure \emph{(see \eqref{app:SRB})}.

\end{theorem}
\bigbreak
In this paper, we illustrate the typical bifurcations from an attracting torus to strange attractors, with a particular example (see Section \ref{s:example1}). 

\section{The example}
\label{s:example1}

We construct an explicit two parametric family of vector fields $f_{(\tau_1, \tau_2)}$ in $\EU^3\subset \RR^4$   whose organizing center satisfies \textbf{(P1)--(P5)}. 
Our construction is based on properties of \emph{differential equations with symmetry} (see  \eqref{lift1}); we  also refer the reader to \cite{GS, GH, GW} for more information on the subject.
\subsection{The system}
\label{s:example}
\medbreak
For $\tau_1, \tau_2 \in \,[0,1]$, our object of study is the two-parameter family of vector fields on $\RR^{4}$ $$x=(x_1,x_2,x_3,x_4)\in\RR^4 \quad \mapsto \quad  f_{(\tau_1, \, \tau_2)}(x)$$
defined for each $x=(x_1,x_2,x_3,x_4)\in \RR^4$ by
\begin{equation}\label{example}
\left\{
\begin{array}{l}
\dot{x}_{1}=x_{1}(1-r^2)-\omega x_2-\alpha x_1x_4+\beta x_1x_4^2 +\textcolor{blue}{\tau_2 x_1x_3x_4} \\ 
\dot{x}_{2}=x_{2}(1-r^2)+{\omega}x_1-\alpha x_2x_4 +\beta x_2x_4^2 \\ 
\dot{x}_{3}=x_{3}(1-r^2)+\alpha x_3x_4+\beta x_3x_4^2+\textcolor{magenta}{\tau_1 x_4^3} -\textcolor{blue}{\tau_2  x_1^2x_4}\\ 
\dot{x}_{4}=x_{4}(1-r^2)-\alpha (x_3^2-x_1^2-x_2^2)-\beta x_4(x_1^2+x_2^2+x_3^2)-\textcolor{magenta}{\tau_1 x_3x_4^2}\\
\end{array}
\right.
\end{equation}
where

$$\dot{x}_i=\frac{\partial x_i}{\partial t}, \qquad r^2=x_{1}^{2}+x_{2}^{2}+x_{3}^{2}+x_{4}^{2}$$
and

$$
\omega>0, \qquad \beta <0<\alpha, \qquad \beta^2<8 \alpha^2 \qquad \text{and} \qquad |\beta|<|\alpha|.
$$

\begin{remark}
The nature of the perturbations, which depend on $\tau_1$ and $\tau_2$,  has been listed in Appendix B of \cite{RodLab}.
\end{remark}

The unit sphere $\EU^3\subset\RR^{4}$ is invariant under the corresponding flow and every trajectory with nonzero initial condition is forward asymptotic to it (cf. \cite{RodLab}). Indeed, if $\left\langle .\,, . \right\rangle $ denotes the usual inner product in $\RR^4$, then it is easy to check that:

\begin{lemma}
For every $x \in \EU^3$ and  $\tau_1, \tau_2 \in [0,1]$, we have $\left\langle f_{(\tau_1,\tau_2)}(x), x\right\rangle =0$.
\end{lemma}
We are interested in  dynamics on a compact boundaryless manifold, in order to have control of the long-time existence and behaviour of solutions.  
Moreover, the origin is repelling since all eigenvalues of $Df_{(\tau_1,\tau_2)}$ at the origin have positive real part. 
\medbreak
\subsection{The organizing center ($\tau_1=\tau_2=0$)}
The vector field $f_{(0,0)}$ is equivariant under the action of the compact Lie group $\mathbb{SO}(2)(\gamma_\psi)\oplus \ZZ_2(\gamma_2)$, where $\mathbb{SO}(2)(\gamma_\psi)$ and $\ZZ_2(\gamma_2)$ act on $\RR^4$ as
$$\gamma_\psi(x_1, x_2,x_3,x_4)=(x_1\cos \psi -x_2 \sin \psi, x_1\sin \psi +x_2\cos \psi, x_3,x_4), \quad \psi \in [0, 2\pi] $$
given by a phase shift $\theta \mapsto \theta+ \psi$ in the first two coordinates, and 
$$ \gamma_2(x_1, x_2,x_3,x_4)=(x_1, x_2,-x_3,x_4).$$
By construction, $\tau_1$ is the controlling parameter of the $\ZZ_2(\gamma_2)-$symmetry breaking and  $\tau_2$ controls the $\mathbb{SO}(2)(\gamma_\psi)-$symmetry breaking but keeping the $\mathbb{SO}(2)(\gamma_\pi)$--symmetry\footnote{Observe that $\mathbb{SO}(2)(\gamma_\pi) \cong \ZZ_2(\gamma_\pi)$, \emph{i.e.} the action of both Lie groups on $\RR^4$ are isomorphic.}, where
$
\gamma_\pi (x_1, x_2, x_3, x_4)=(-x_1, -x_2, x_3, x_4). $
See the following table for the symmetries preserved according to parameters:
\bigbreak
\bigbreak
\begin{center}
\begin{tabular}{lr}
\hline \hline & \\
Parameters \qquad \qquad  & Symmetries preserved   \\ & \\ \hline \hline
&\\
$\tau_{1}=\tau_2=0$ & $\mathbb{SO}(2)(\gamma_\psi) \oplus  \ZZ_2(\gamma_2) $  \\  & \\ \hline & \\
$\tau_{1}>0\,  \text{ and}\, \tau_2=0$ & $\mathbb{SO}(2)(\gamma_\psi) $   \\  & \\ \hline & \\
$\tau_{1}=0\,  \text{ and}\, \tau_2>0$ & $ \ZZ_2(\gamma_\pi) \oplus  \ZZ_2(\gamma_2) $   \\  & \\ \hline & \\
$\tau_{1}>0\,  \text{ and}\, \tau_1 \gg \tau_2$ \qquad & $\ZZ_2(\gamma_\pi) $   \\  & \\ \hline \hline
\end{tabular}
\end{center}
\begin{center}
\small{Table 1: Types of symmetry-breaking according to the parameters.}
\end{center}

\bigbreak

\begin{figure}[h]
\begin{center}
\includegraphics[height=4.5cm]{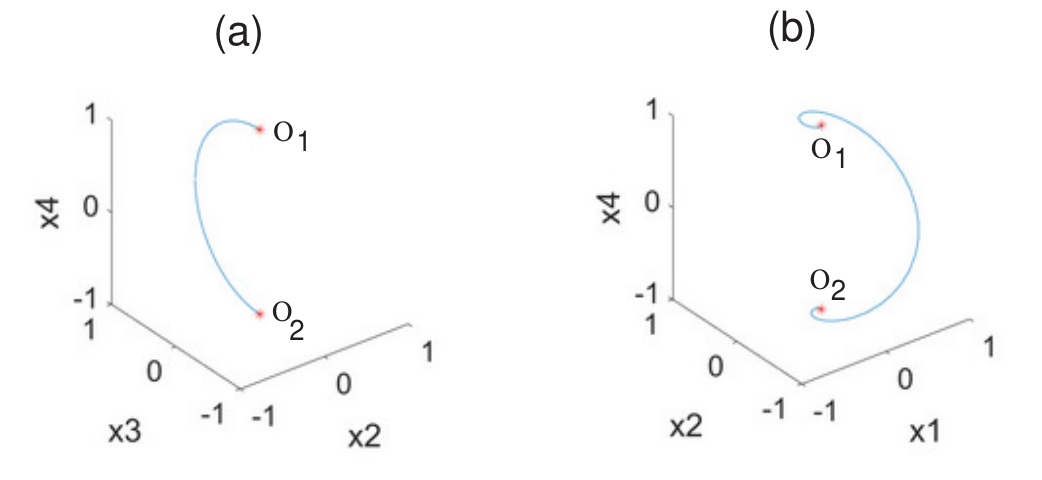}
\end{center}
\caption{\small Organizing center. (a) One 1D-connection from $O_1$ to $O_2$. (b) One trajectory within the 2D-connection from $O_2$ to $O_1$. Flow of \eqref{example} with $\tau_1=\tau_2=0$, $\omega=1$, $\alpha=1$ and $\beta=-0.1$, $t\in [0, 10000]$. (a) Initial condition $(0; 0; 0.01; 0.99)$ near $W^s(O_1)$; projection into the coordinates $(x_2, x_3, x_4)$. (b) Initial condition $(0.1; 0.1; 0; -0.99)$ near $W^u(O_2)$; projection into the coordinates $(x_1, x_2, x_4)$.  } 
\label{orientation_exemplo1}
\end{figure}

\bigbreak
When restricted to the sphere $\EU^3$, for every $\tau_1, \tau_2 \in [0,1]$, the flow of $f_{(\tau_1, \tau_2)}$ has 
 two equilibria 
$$O_1 =(0,0,0,+1) \quad \quad \text{and} \quad \quad O_2 = (0,0,0,-1), $$
which are hyperbolic saddle-foci with different \emph{Morse indices} (dimension of the unstable manifold). 
The linearization of $f_{(0,0)}$ at $O_1$ and $O_2$ has eigenvalues
$$ -(\alpha-\beta) \pm \omega i, \,\,  \alpha+\beta \qquad \text{and} \qquad (\alpha + \beta)\pm \omega i, \,\,,  -(\alpha-\beta)$$
respectively. As depicted in Figure \ref{orientation_exemplo1}, when restricted to $\EU^3$, in these coordinates, the 1D-connections are given by: 
 $$W^u(O_1) \cap \EU^3= W^s(O_2) \cap \EU^3=\text{Fix}(\mathbb{SO}(2)(\gamma_\psi))\cap \EU^3= \{(x_1,x_2,x_3,x_4): x_1=x_2=0, x_3^2 + x_4^2 = 1\}$$
and the 2D-connection is contained in
$$W^u(O_2) \cap \EU^3= W^s(O_1) \cap \EU^3=\text{Fix}(\ZZ_2(\gamma_2))\cap \EU^3= \{(x_1,x_2,x_3,x_4): x_1^2 + x_2^2 + x_4^2 = 1, x_3=0\}.$$

The two-dimensional invariant manifolds are contained in the two-sphere $\text{Fix}(\ZZ_2(\gamma_2 ))\,\cap\, \EU^3.$ It is precisely this symmetry that forces the two-invariant manifolds $W^u(O_2)$ and $W^s(O_1)$ to coincide. In what follows, we denote by $\Gamma$ the \emph{heteroclinic network} formed by the two equilibria, the two connections $[O_1 \rightarrow O_2]$ and the sphere $[O_2 \rightarrow O_1]$ (see Figure~\ref{orientation_exemplo1}). The network may be decomposed into two cycles.

\medbreak

\begin{figure}[ht]
\begin{center}
\includegraphics[height=4.4cm]{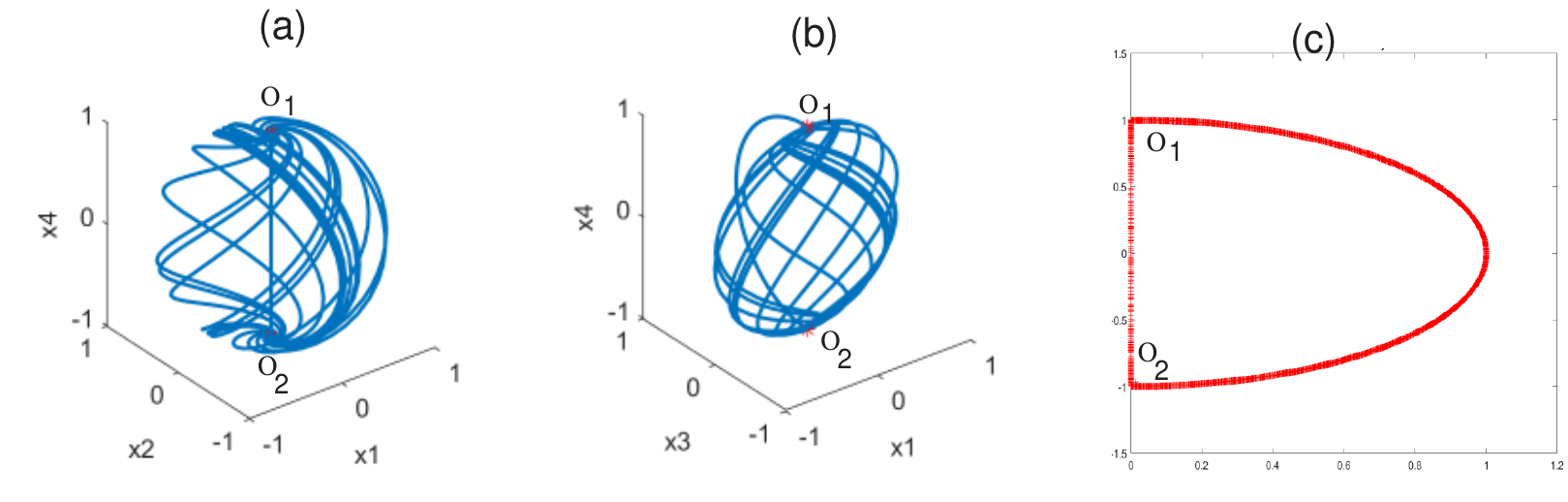}
\end{center}
\caption{\small Flow of \eqref{example} of the trajectory with initial condition $(0.1; 0.1; 0; -0.99)$ near $W^u(O_2)$, with $\tau_1=\tau_2=0$, $\omega=1$, $\alpha=1$ and $\beta=-0.1$, $t\in [0, 10000]$. (a) Projection into the coordinates $(x_1, x_2, x_4)$. (b) Projection into the coordinates $(x_1, x_3, x_4)$. (c) Projection into the section $x_1=x_2=0$.} 
\label{asymptotically_stable}
\end{figure}

  Proposition 1 of \cite{RodLab} shows that, keeping $\tau_1=\tau_2=0$, the equilibria $O_1$ and $O_2$ have the same \emph{chirality} (details in \eqref{chirality_def}). Therefore:
\begin{lemma}
If $\tau_1=\tau_2=0$,  the flow of \eqref{example} satisfies \textbf{(P1)--(P5)} described in Section \ref{s:setting}. 
\end{lemma}
In summary, when $\tau_1=\tau_2=0$, the flow of \eqref{example} exhibits an asymptotically stable heteroclinic network $\Gamma$ associated to $O_1$ and $O_2$, numerically shown in Figure \ref{asymptotically_stable}. Throughout the construction and discussion, the parameters $\tau_1$ and $\tau_2$   play the role of $A$ and $\lambda$, respectively, of \textbf{(P6)--(P7)}, after possible rescaling.

\subsection{Notation} As in Subsection \ref{notation1}, we set the following notation:
\begin{equation}
\label{constants}
C_1 = C_2 = \alpha - \beta>0, \qquad E_1 = E_2 = \alpha + \beta>0, \qquad \delta_1= \delta_2 = \frac{\alpha - \beta}{\alpha + \beta }>1
 \end{equation}
and
\begin{equation}
\label{constants2}
K = \frac{2 \alpha}{(\alpha +\beta)^2}>0  \qquad \text{and} \qquad  K_\omega= \frac{2 \alpha \omega  }{(\alpha +\beta)^2}>0.
 \end{equation}

\subsection{$\ZZ_2(\gamma_2)$--symmetry breaking ($\tau_1>\tau_2=0$)}
In this scenario,  the heteroclinic network $\Gamma$ is broken because the symmetry $\ZZ_2(\gamma_2)$ is broken.
The flow of $f_{(\tau_1, 0)}$ leaves the unit sphere $\EU^3$ invariant and globally attracting since the perturbations are tangent to  $\EU^3$ \cite[Appendix B]{RodLab}. 
We are going to present analytical evidences that  an attracting two-torus is born.

Since the system $\dot{x}=f_{(\tau_1, 0)}(x)$ is still $\mathbb{SO}(2)(\gamma_\psi)$--equivariant (see Table 1), we may define the quotient flow on $\EU^3 / \mathbb{SO}(2)(\gamma_\psi)$ (see \cite{Schwarz}) and we get the following differential equation:
\begin{equation}\label{exampleR3}
\left\{
\begin{array}{l}
\dot{\rho}=\rho(1-R^2)-\alpha \rho \, x_4+\beta   \rho \,x_4^2 \\ \\
\dot{x}_{3}=x_{3}(1-R^2)+ \alpha x_3x_4+ \beta x_3 x_4^2+\tau_1 x_4^3 \\\\
\dot{x}_{4}=x_{4}(1-R^2)-\alpha (x_3^2-\rho^2)-\beta x_4(\rho^2+x_3^2)-\tau_1 x_3x_4^2\\
\end{array}
\right.
\end{equation}
where
$$R^2=\rho^{2}+x_{3}^{2}+x_{4}^{2} \qquad \text{and} \qquad \rho^2=x_1^2+x_2^2.$$
The equations \eqref{exampleR3}, restricted to the unit two-sphere $\EU^2$ (\emph{i.e.} $R^2=1$), simplify to
\begin{equation}\label{exampleR3.2}
\left\{
\begin{array}{l}
\dot{\rho}=\alpha \rho \, x_4+\beta   \rho \,x_4^2 \\ \\
\dot{x}_{3}= \alpha x_3x_4+ \beta x_3 x_4^2+\tau_1 x_4^3 \\\\
\dot{x}_{4}=-\alpha (x_3^2-\rho^2)-\beta x_4(\rho^2+x_3^2)-\tau_1 x_3x_4^2,\\

\end{array}
\right.
\end{equation}
a differential equation  that may be \emph{reduced} to the following planar system:
\begin{equation}
\label{exampleR3.3}
\left\{
\begin{array}{l}
\dot{x}_{3}=\alpha x_3x_4 + \beta x_3 x_4^2+\tau_1 x_4^3 \\ \\
\dot{x}_{4}=\alpha (1-2x_3^2-x_4^2) + \beta x_4(x_4^2-1)-\tau_1 x_3x_4^2.\\
\end{array}
\right.
\end{equation}

\begin{figure}[h]
\begin{center}
\includegraphics[height=6.0cm]{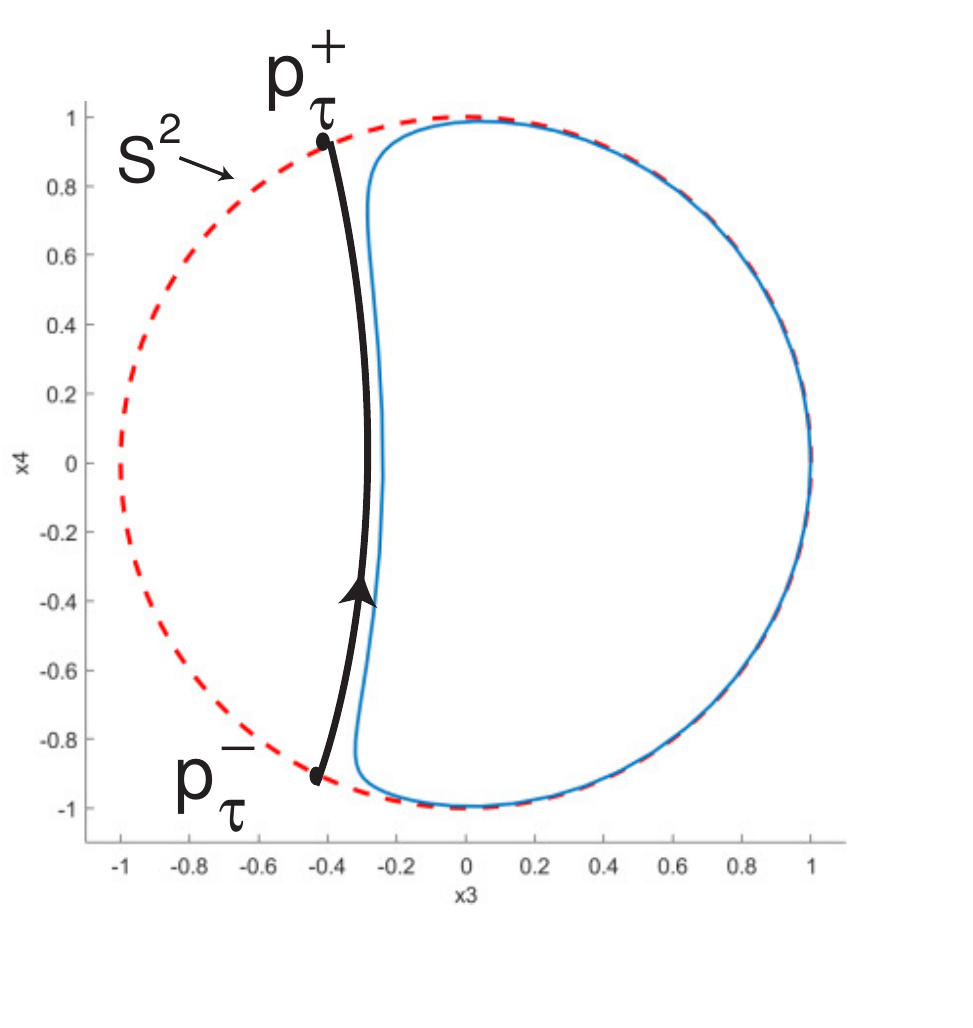}
\end{center}
\caption{\small Flow of \eqref{exampleR3.3} of the trajectory with initial condition $(x_3, x_4)=(0; -0.99)$ near $W^u(p^-_\tau)$, with $\alpha=1$, $\beta=-0.1$ and $\tau_1=0.5$, $t \in [0, 10000]$. The blue closed curve of \eqref{exampleR3.3} is stable. The red dashed line represents the unit circle. } 
\label{torus_dstool}
\end{figure}

\begin{figure}[h]
\begin{center}
\includegraphics[height=5.0cm]{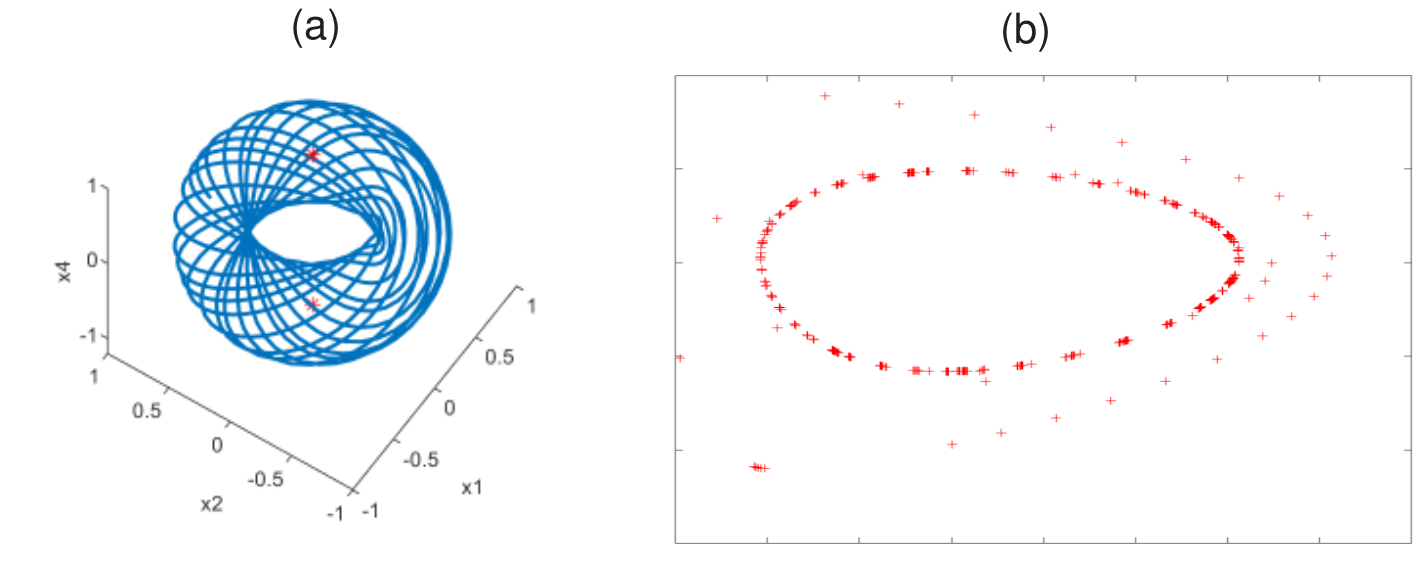}
\end{center}
\caption{\small Projection of the flow of \eqref{example} of the trajectory with initial condition $(0.1; 0.1; 0; -0.99)$ near $W^u(O_2)$, with $\tau_1=0.5$ $\tau_2=0$, $\omega=1$, $\alpha=1$ and $\beta=-0.1$, $t\in [0, 10000]$. (a) Projection into the coordinates $(x_1, x_2, x_4)$. (b) Projection into the section $x_1=x_2=0$. } 
\label{torus_dstool2}
\end{figure}

For $\tau_1=0$, the points $p^\pm_0=(0, \pm 1)\in \EU^1$ are hyperbolic equilibria for \eqref{exampleR3.3}. For $\tau_1\neq 0$, let  $p^\pm_\tau$ be  their \emph{hyperbolic continuation}.  Using the Poincar\'e-Bendixson Theorem, Aguiar \cite{Aguiar_tese} proved that:
\begin{lemma}[\cite{Aguiar_tese}, adapted]
\label{sinks_lemma}
For $\tau_1>0$, the flow of system \eqref{exampleR3.3} has one stable periodic solution, which emerges from the breaking of the attracting network associated to the equilibria $p^\pm_0$. The unstable manifold of $p^-_\tau$ does not intersect the stable manifold of $p^+_\tau$.
\end{lemma}
The  stable periodic solution of Lemma \ref{sinks_lemma} is illustrated in Figure \ref{torus_dstool}.
By the $\mathbb{SO}(2)(\gamma_\psi)$--equivariance, this sink lifts to an attracting torus (details of the \emph{lifting process} is given in  \eqref{lift1}).
Therefore:
\begin{corollary}
\label{torus_cor}
For $\tau_1>0$ and $\tau_2=0$, close to the ``ghost'' of the attracting network $\Gamma$, the flow of \eqref{example}  has an attracting invariant two-torus, which is normally hyperbolic. 
\end{corollary}
  From the theory for normally hyperbolic manifolds developed in \cite{HPS}, the torus persists under small smooth perturbations.  The numerical evidence of Figure \ref{torus_dstool2}(a) suggests that the dynamics restricted to the torus is \emph{quasi-periodic}.

\subsection{$\ZZ_2(\gamma_2)$ and $\mathbb{SO}(2)(\gamma_\psi)$--symmetry breaking ($\tau_1\gg \tau_2>0$)}

We now  explore the case $\tau_1\gg \tau_2>0$.
Once again, the flow of $f_{(\tau_1, \tau_2)}$ leaves the unit sphere $\EU^3$ invariant and globally attracting since the perturbations are tangent to  $\EU^3$ -- \cite[Appendix B]{RodLab}. 
Although we break the $\mathbb{SO}(2)(\gamma_\psi)$--equivariance, the $\mathbb{SO}(2)(\gamma_\pi)$--symmetry is preserved. This is why the connections lying in  $x_1=x_2=0$ persist ($\Rightarrow$ \textbf{(P6)} holds).
Using now Lemma \ref{sinks_lemma}, by construction, we have:
\begin{lemma}
For $\tau_1\geq 0$ and $\tau_2 > 0$ small enough such that $\tau_1\gg \tau_2>0$, the flow of \eqref{example} satisfies \textbf{(P6)--(P7)}. 
\end{lemma}

\begin{figure}[ht]
\begin{center}
\includegraphics[height=15.0cm]{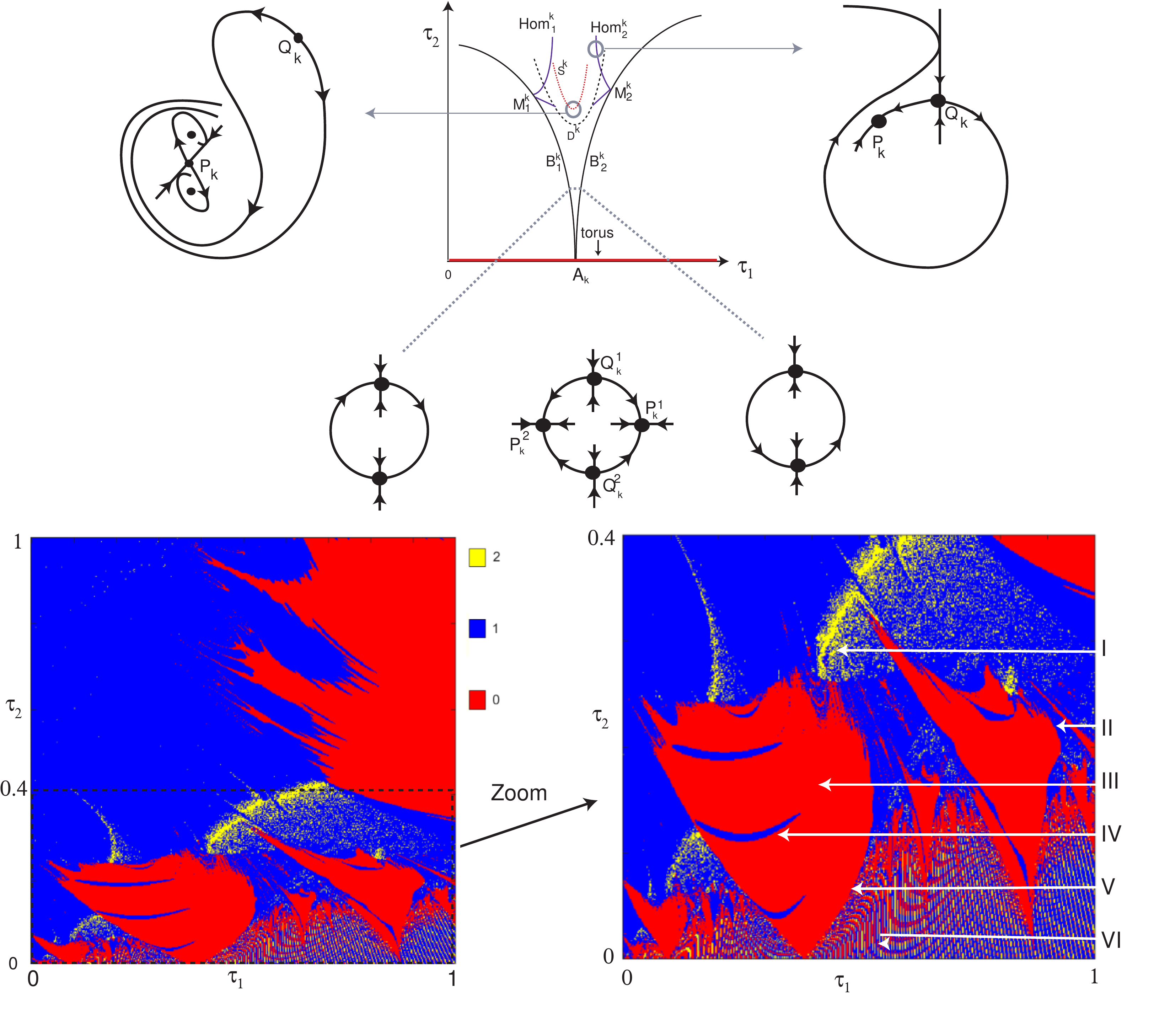}
\end{center}
\caption{\small Bifurcation diagram for equation \eqref{example} with $\alpha=1$, $\beta=-0.1$ and $\omega=1$, corresponding to the trajectory with initial condition $(0.1; 0.1; 0; -0.99)$ near $W^u(O_2)$, $t\in [0,3750]$. Upper part: theoretical scheme \cite{Aronson}. Lower part: (real) bifurcation diagram.  
Each point is colored according to the number of non-negative Lyapunov exponents of the orbit: yellow (2), blue (1) and red (0), corresponding to an $\omega$--limit including a strange attractor or an attracting two-torus, a limit cycle or a fixed point, respectively. } 
\label{DiagBif1}
\end{figure}

 Numerical simulations of  \eqref{example} for $\tau_1\gg \tau_2>0$ suggest the existence of regular and  chaotic behaviour in the region of transition from regular dynamics (attracting torus) to rotational saturated horseshoes (see \eqref{ss: saturated horseshoe} and \eqref{ss: rotational horseshoe}). Chaotic attractors with  one positive Lyapunov exponent seem to exist, as suggested by the yellow regions of Figure \ref{DiagBif1}.
Using the Matlab software (R2019b, MathWorks, Natick, MA, USA), we have been able to compute the bifurcation diagram presented in Figure \ref{DiagBif1}, whose analysis is the goal of next section.

\section{Torus-breakdown and strange attractors: theory and numerics}
\label{interpretation}

The destruction of the torus of Corollary \ref{torus_cor} takes place according to the scenario described in  \cite{AS91, AHL2001, AH2002, Aronson}.  In this section, we describe the (generic) mechanism to break an attracting two-dimensional torus and we relate this theory with our numerics, which agree quite well.
Before going further, we introduce some terminology. 
\medbreak

Let $\mathcal{T}$ be a neighborhood of the Bykov attractor $\Gamma$ which exists for \eqref{example} when $\tau_1=\tau_2=0$.   Let $\Sigma$ be a cross section to $\Gamma$. There is $\varepsilon>0$ (small) such that the first return map $\mathcal{F}_{(\tau_1, \tau_2)} $ to a subset of ${\Sigma}$ is well defined, for $\tau_1, \tau_2<\varepsilon$ (cf. \cite{Rodrigues2019}). We assume that the intersection of the torus with the cross section $\Sigma$ is a curve diffeomorphic to a circle, as depicted in Figure \ref{torus_dstool2}(b).

\medbreak

The choice of parameters in Section \ref{s:example1} lets us build  the bifurcation diagram in the plane of the parameters $(\tau_1, \tau_2)$ in the domain $$\{0 \leq \tau_1 < \varepsilon, \quad 0 < \tau_2 < \tau_2^0\},$$ for $0<\tau_2^0 \ll \varepsilon$.  Within this region, it is possible to define an \emph{Arnold tongue} \cite{Aronson}, denoted by $\mathcal{T}_k$, adjoining the horizontal axis at a point $A_k$ where $k$ is an  integer. Inside this \emph{tongue} (\emph{resonant wedge}), for small $\tau_2$, there coexist at least a pair of fixed points for the Poincar\'e map $\mathcal{F}_{(\tau_1, \tau_2)} $, whose corresponding trajectories share the same \emph{rotation number} \cite{AS91, Herman}.  As illustrated in Figure \ref{DiagBif1} (upper part), we suppose the existence of two pairs of fixed points: $Q_k^1$, $Q_k^2$ (saddles) and $P_k^1$, $P_k^2$ (sinks).  
\medbreak
The borders of $\mathcal{T}_k$ are  bifurcation curves $B_1^k$ and $B_2^k$ on which each pair of fixed points merge into a \emph{saddle-node}. These curves might touch the corresponding curves of other tongue, meaning that there are parameter values for which  periodic solutions with different rotation number might coexist.
  The points $M_1^k$ and $M_2^k$ correspond to \emph{homoclinic cycles to a saddle-node}: below these points, in $B_1^k$ and $B_2^k$,  the limit set of $W^u(Q_k^1)$ is the saddle-node itself. Above the points $M_1^k$ and $M_2^k$, the maximal invariant set is not homeomorphic to a circle.  In the bifurcation diagram, there is also a curve, say $D^k$, above which the invariant torus no longer exists due to a \emph{period doubling bifurcation} \cite{Anishchenko}.   After the period doubling has occurred, the torus is destroyed.
 
 \medbreak

Continuing the process of dissecting an Arnold tongue, the authors of \cite{AS91, Anishchenko} describe generic mechanisms  by which the invariant and attracting torus is destroyed. Two of them are revived in the next result and involve  homoclinic tangencies -- routes [PA] and [PB] of \cite{Anishchenko}. 

\begin{theorem}[\cite{AS91, Anishchenko}, adapted]
\label{three_lines}
For $K_\omega^0>0$ fixed, in the bifurcation diagram $\dpt \left(\tau_1, {\tau_2}\right)$, within $\mathcal{T}_k$, 
\begin{itemize}  
\item[(1)]   there are two curves $\text{Hom}_1^k$ and $\text{Hom}_2^k$ corresponding to a homoclinic tangency associated to a dissipative periodic point of the first return map $\mathcal{F}_{(\tau_1, \tau_2)}$.
\item[(2)] there is one curve $S^k$ corresponding to a homoclinic tangency (of third class) associated to a dissipative periodic point of the first return map $\mathcal{F}_{(\tau_1, \tau_2)}$.
\end{itemize}
 \end{theorem}

The lines $S^k$, $\text{Hom}_1^k$, $\text{Hom}_2^k$ are shown in Figure \ref{DiagBif1}. \medbreak

\begin{figure}[ht]
\begin{center}
\includegraphics[height=7.5cm]{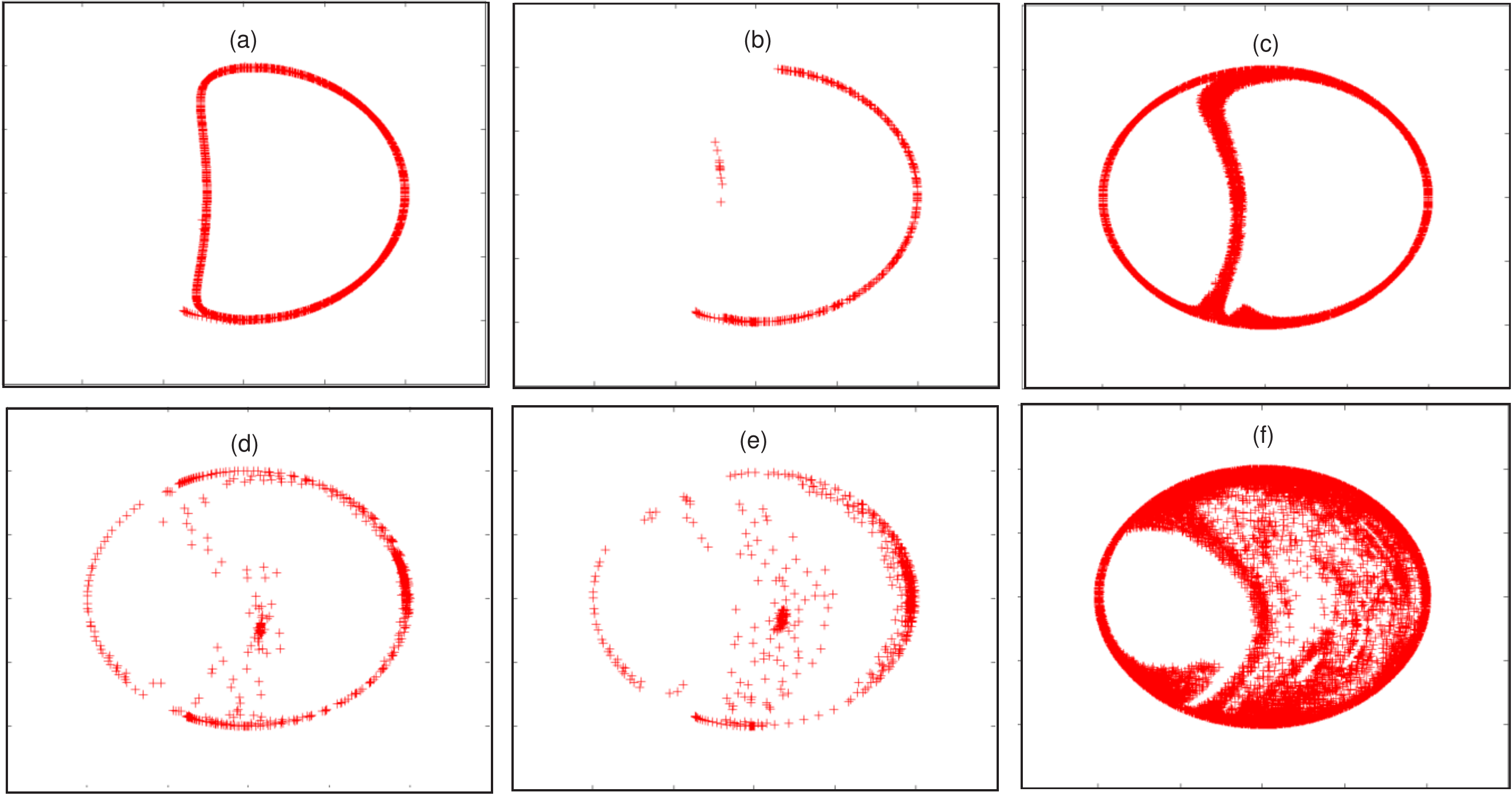}
\end{center}
\caption{\small Projection of the flow of \eqref{example} of the trajectory with initial condition $(0.1; 0.1; 0; -0.99)$ with $\alpha=1$, $\beta=-0.1$,   $\tau_1=0.3$, $\omega=1$ and $\tau_2=0$ (a);  $\tau_2=0.1$ (b),  $\tau_2=0.2$ (c), $\tau_2=0.3$ (d),  $\tau_2=0.4$ (e),  and $\tau_2=0.5$ (f).} 
\label{strange_attractor1}
\end{figure}

Fix $i\in \{1,2\}$ and assume that $Q_k^i\equiv Q_k$.  Along the bifurcation curves $\text{Hom}_1^k$ and $\text{Hom}_2^k$, one observes a homoclinic contact of  the components $W^s(Q_k)$ and $W^u(Q_k)$, where $Q_k$ is a dissipative saddle\footnote{One saddle $O$ is dissipative if $0<|\det D\mathcal{F}_{(\tau_1, \tau_2)}(O)|<1$ (\emph{i.e.}  $\mathcal{F}_{(\tau_1, \tau_2)}$ is contracting for any small neighbourhood of $O$. Note that, for small $\tau_1, \tau_2>0$, the first return map $\mathcal{F}_{(\tau_1, \tau_2)}$ is contracting. }.
The curves $\text{Hom}_1^k$ and $\text{Hom}_2^k$ divide the region above  $D^k$ into two regions with simple and complex dynamics. In the zone above the curves $\text{Hom}_1^k$ and $\text{Hom}_2^k$, there is a fixed point $Q_k$ exhibiting a transverse homoclinic intersection, and thus the corresponding  map $\mathcal{F}_{(\tau_1, \tau_2)}$ exhibits nontrivial hyperbolic chaotic sets (horseshoes). Other stable points of large period exist in the region above the curves $\text{Hom}_1^k$ and $\text{Hom}_2^k$ since the homoclinic tangencies arising in these lines are generic -- \emph{Newhouse phenomena} \cite{Colli98, Newhouse79}. 
Using now \cite{MV93}, there exists a positive measure set  $\Delta$ of parameter values, so that for every $\tau_1/\tau_2\in \Delta$, the map $\mathcal{F}_{(\tau_1, \tau_2)}$ admits a strange attractor  of H\'enon-type with an ergodic SRB measure (cf. \eqref{app:SRB}).  This is observable in the two yellow regions leaving the \emph{Arnold tongues} of Figure \ref{DiagBif1} (lower part). The formation of the H\'enon-like strange attractor is suggested in Figure \ref{strange_attractor1}.

The curve $S^k$ of Theorem \ref{three_lines} corresponds to a homoclinic tangency of third class, meaning that there are tangencies associated to the fixed points  emerged from the period-doubling bifurcation at $D^k$.  This line corresponds to the boundary between the red and blue lines in Figure \ref{DiagBif1} (lower image), within the resonant wedge.  We conjecture that this line is a consequence of an exponentially small wedge associated to a \emph{Bogdanov-Takens} bifurcation \cite[Fig. 3]{Yagasaki}. After crossing  the curve $D_k$ (from below), the invariant curve no longer exists. Of course, the individual solutions (orbits) remain smooth.

\medbreak

Yellow narrow regions of Figure \ref{DiagBif1}, near the horizontal axis, are due to the existence of a set of parameters for which the torus-flow  is \emph{irrational}; the corresponding orbit is unlocked and winds without bound around the torus. In Figure \ref{DiagBif1}, we also observe \emph{frequency locking regions} (red regions), regions dominated by a sink (blue) as well as regions with positive entropy (yellow). The number of connected components with which the strange attractors intersect the cross-section is not specified nor is the size of their basins of attraction.

\begin{figure}[h]
\begin{center}
\includegraphics[height=4.8cm]{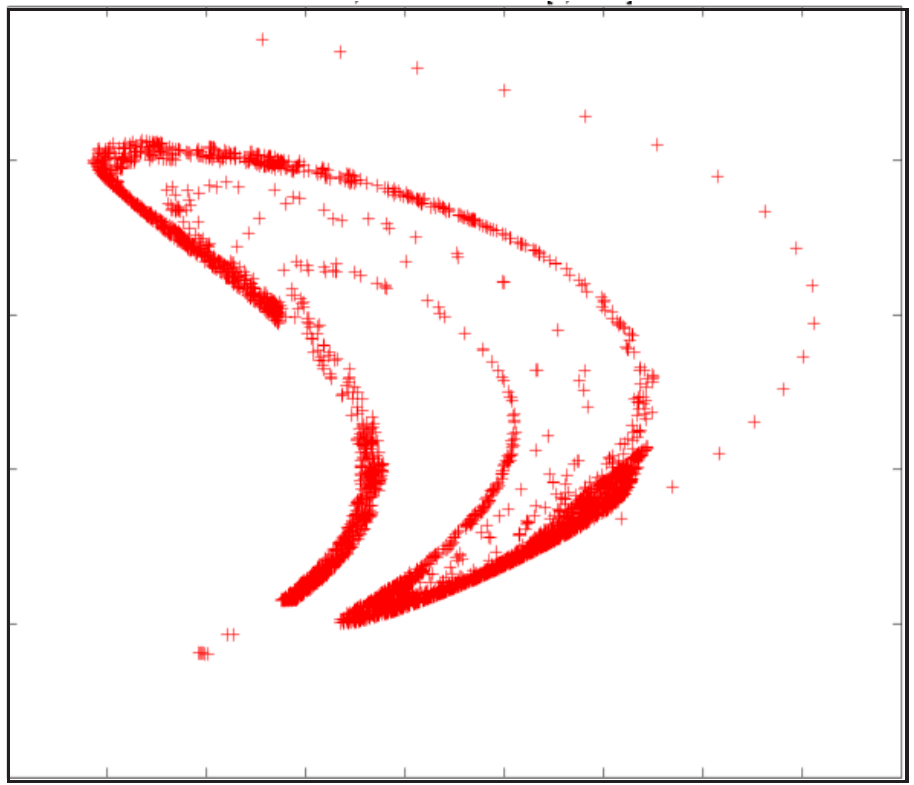}
\end{center}
\caption{\small Breaking the invariant curve of Figure \ref{torus_dstool2}(b) on the section $x_1=x_2=0$. Projection of the flow of \eqref{breaking_all} of the trajectory with initial condition $(0.1; 0.1; 0; -0.99)$ and $t\in [0, 10000]$. The perturbation breaks all the symmetries, as well as the sphere invariance. } 
\label{strange_attractor2}
\end{figure}

For $\tau_1, \tau_2 >0$, the flow associated to $f_{(\tau_1,\tau_2)}$ leaves $\EU^3$ invariant and globally attracting, which explains the regularity of borders of the different images in Figure \ref{strange_attractor1}. In Figure \ref{strange_attractor2}, we give an additional generic image of the typical picture of \emph{H\'enon-like attractor} which appears near the Bykov attractor when all symmetries are broken, as well as the  sphere-invariance. The plot has been performed for the vector field $f_{(0,0)}$ of \eqref{example} with the following perturbing term:
\begin{equation}
\label{breaking_all}
(x_1x_3x_4, \,\,   -x_1 x_2^2,\, \,   x_3^3,\,\,  - x_1x_3x_4).
\end{equation}
The resulting vector field breaks all symmetries, all well as the sphere invariance. The unstable manifold of the saddle $Q_k$ has crossed the non-leading stable manifold of the periodic orbit.

\subsection*{Possible interpretation of Figure \ref{DiagBif1}}

 \begin{eqnarray*}
\textbf{I} \quad &\to& \quad \text{Homoclinic bifurcations; H\'enon-like strange attractors}.\\
\textbf{II}\quad &\to& \quad \text{Sink}.\\
\textbf{III}  \quad &\to& \quad \text{Resonant wedge (Arnold tongue)}.\\
\textbf{IV} \quad &\to& \quad \text{Hopf bifurcation}.\\
\textbf{V}\quad &\to& \quad \text{Saddle-node bifurcation (border of the Arnold tongue)}.\\
\textbf{VI}  \quad &\to& \quad \text{Irrational torus (thin yellow region)}.\\
\end{eqnarray*}
 
\subsection*{Technicalities of the numerics}

Since $W^u(O_2)$ plays an essential role in the construction of the H\'enon-like strange attractor \cite{Rodrigues2019}, we chose  $(0.1; 0.1; 0; -0.99) \in \RR^4$ to grasp the main dynamical properties of the maximal attracting set of \eqref{example}. 
Although the system \eqref{example} lives in $\RR^4$, the analysis may be performed in the sphere $\EU^3$ since it is globally attracting.  According to \cite[pp. 287]{Wolf85}, for a three-dimensional continuous dissipative flow, the only possible spectra and the attractors they describe depend of the sign of their Lyapunov exponents:
\begin{eqnarray*}
{(-,\, -,\, -)} \quad &\to& \quad \text{the $\omega-$limit of the corresponding orbit contains a fixed point};\\
{(0,\, -,\, -)} \quad &\to& \quad \text{the $\omega-$limit of the corresponding orbit is a limit cycle};\\
{(+,\, 0,\, -)} \quad &\to& \quad \text{the $\omega-$limit of the corresponding orbit is a chaotic attractor};\\
{(0,\, 0,\, -)} \quad &\to& \quad \text{the $\omega-$limit of the corresponding orbit is an attracting 2-torus}.\\
\end{eqnarray*}

 The parameter plane $(\tau_1, \tau_2)$ of Figure \ref{DiagBif1} is scanned with a sufficiently small step along each coordinate axes. The software evaluates at each parameter value how many Lyapunov exponents along the orbit with initial condition $(0.1; 0.1; 0; -0.99)$ are non-negative. Then, the parameter is painted according to the following rules: red for $0$, blue for $1$, yellow for $2$.  To estimate the complete Lyapunov spectra, we use the algorithm for differential
equations introduced in \cite{Wolf85} with a Taylor series integrator.

\begin{remark}
Observe that yellow regions of Figure \ref{DiagBif1} mean that the $\omega-$limit of the corresponding orbit may be either a strange attractor (I) or an attracting two-torus (VI).  The first seems to be prevalent; the second is not. 
\end{remark}


\section{Discussion and final remarks}
\label{s:discussion}

The goal of this paper is to construct explicitly a two-parameter family of polynomial differential equations $\dot{x}=f_{(\tau_1, \tau_2)}(x)$ in the three-dimensional sphere $\EU^3$, in which each parameter controls a type of symmetry-breaking. Depending on the parameters, different dynamical regimes have been identified both analytically and numerically. We have stressed the emergence of strange attractors from an attracting heteroclinic network, a  \emph{route  to chaos} which has been a recurrent concern on nonlinear dynamics during the last decades. 
 Along this discussion we compare our results to what is known for other models in the literature.

The flow of $\dot{x}=f_{(0,0)}(x)$ has an attracting heteroclinic network $\Gamma$ with a non-empty basin of attraction $\mathcal{U}$.
 We have studied  the global transition of the dynamics from $\dot{x}=f_{(0,0)}(x)$ to a  smooth two-parameter family $\dot{x}=f_{(\tau_1,\tau_2)}(x)$ that  breaks part of the network. 
For small perturbations, the set $\mathcal{U}$ is still positively invariant. When $\tau_1, \tau_2 \neq 0$,
the one-dimensional connections persist due to the remaining symmetry and, as  a consequence of Kupka-Smale Theorem,  the 
two-dimensional invariant manifolds are generically transverse (either intersecting or not).

When $\tau_2>\tau_1 \geq 0$, the two-dimensional invariant manifolds intersect transversely, giving rise to a complex network, that consists of a union of Bykov cycles \cite{Bykov00}, contained 
in $\mathcal{U}$. The dynamics in the maximal invariant set contained in $\mathcal{U}$, contains, but does not coincide with, the suspension of horseshoes accumulating on the heteroclinic network described in \cite{ACL05, KLW, LR, Rodrigues2, Rodrigues3, RodLab}. In addition, close to the organizing center ($\tau_1=\tau_2=0$), it contains infinitely many heteroclinic tangencies and attracting limit cycles with long  periods, coexisting with sets with positive entropy, giving rise the so called \emph{quasi-stochastic attractors}.  The sinks in a quasi-stochastic attractor have long periods and narrow basins of attraction, and they are hard to be observed in applied problems  \cite{Afraimovich83, Gonchenko97}. 


The scenario $\tau_1>\tau_2\geq 0$  corresponds to the case where the two-dimensional invariant manifolds do not intersect.  Although the network associated to the equilibria is destroyed,  complex dynamics appears near its ``ghost''. 
In the present article, it is shown that the perturbed system may manifest regular behaviour corresponding to the existence of a smooth invariant torus, and may also have chaotic regimes. In the region of transition from regular behaviour to chaotic dynamics, using known results about \emph{Arnold tongues}, we illustrate the existence of lines with homoclinic tangencies to dissipative periodic solutions, responsible for the existence of persistent strange attractors nearby (cf. region (I) of Figure \ref{DiagBif1}).  Numerics agree quite well with the theory described in \cite{Rodrigues2019}.
Persistence of chaotic dynamics is physically relevant because it means that the phenomenon is numerically observable with positive probability.

In the fully asymmetric case, the  general study of \eqref{example} seems to be analitically untreatable. We have been able to predict qualitative features of the dynamics of the perturbed vector field by assuming that the perturbation is very close to the organizing center. 
Symmetry plays two roles:
first, it creates flow-invariant subspaces where non-transverse heteroclinic connections are persistent, and hence cycles are robust in this context;
second, we use the proximity of the fully symmetric case to capture global dynamics.
Symmetry constrains the  geometry of the invariant manifolds of the saddle-foci and allows us some control of their relative positions. This is an important advantage of studying systems close to symmetry.

In a subsequent paper we will treat some other aspects of this class of examples, in particular partial mode-locking, frequently associated with the existence of homotopically non-trivial invariant circles on the torus. We conjecture that the boundaries of partial mode-locked regions involve \emph{Bodganov-Takens bifurcations} in the parameter space $(\tau_1, \tau_2, K_\omega)$.

\section*{Acknowledgements}

AR was partially supported by CMUP (UID/MAT/00144/2019), which is funded by FCT with national (MCTES) and European structural fundsds through the programs FEDER, under the partnership agreement PT2020. AR also acknowledges financial support from Program INVESTIGADOR FCT (IF/00107/2015).

\appendix

\section{Glossary}\label{Definitions}

For $\varepsilon>0$ small enough, consider the two-parameter family of $C^3$--smooth autonomous differential equations
\begin{equation}
\label{general2a}
\dot{x}=f_{(\tau_1, \tau_2)}(x)\qquad x\in \EU^3  \qquad \tau_1, \tau_2 \in [0, \varepsilon] 
\end{equation}
where $\EU^3$ denotes the unit sphere, endowed with the usual topology. Denote by $\varphi_{(\tau_1, \tau_2)}(t,x)$, $t \in \RR$, the associated flow.

\subsection{Attracting set}
\label{attracting set}

A subset $\Omega$ of a topological space $\mathcal{M}$ for which there exists a neighborhood $U \subset \mathcal{M}$ satisfying $\varphi(t,U)\subset U$ for all $t\geq 0$ and $\bigcap_{t\,\in\,\RR^+}\,\varphi(t,U)=\Omega$ is called an \emph{attracting set} by the flow $\varphi$, not necessarily connected. Its basin of attraction, denoted by $\textbf{B}(\Omega)$ is the set of points in $\mathcal{M}$ whose orbits have $\omega-$limit in $\Omega$. We say that $\Omega$ is \emph{asymptotically stable} (or that $\Omega$ is a \emph{global attractor}) if $\textbf{B}(\Omega)=\EU^3\backslash\{{\textbf{O}}\}$. An attracting set is said to be \emph{quasi-stochastic} if it encloses periodic solutions with different Morse indices, structurally unstable cycles, sinks and saddle-type invariant sets.

\subsection{Heteroclinic phenomenon and Bykov cycle}
\label{ss:Bylov_cycle}
Suppose that $O_1$ and $O_2$ are two hyperbolic saddle-foci of \eqref{general2a} with different Morse indices (dimension of the unstable manifold). There is a {\em heteroclinic cycle} associated to $O_1$ and $O_2$ if
$W^{u}(O_1)\cap W^{s}(O_2)\neq \emptyset$ and  $W^{u}(O_2)\cap W^{s}(O_1)\neq \emptyset.$ For $i, j \in \{1,2\}$, the non-empty intersection of $W^{u}(O_i)$ with $W^{s}(O_j)$ is called a \emph{heteroclinic connection} between $O_i$ and $O_j$, and will be denoted by $[O_i \rightarrow  O_j]$. Although heteroclinic cycles involving equilibria are not a generic feature within differential equations, they may be structurally stable within families of systems which are equivariant under the action of a compact Lie group $\mathcal{G}\subset \mathbb{O}(n)$, due to the existence of flow-invariant subspaces \cite{GH}.

A heteroclinic cycle between two hyperbolic saddle-foci of different Morse indices, where one of the connections is transverse (and so stable under small perturbations) while the other is structurally unstable, is called a Bykov cycle. A \emph{Bykov network} is a connected union of heteroclinic cycles, not necessarily in finite number. 
We refer to \cite{HS} for an overview of heteroclinic bifurcations and substantial information on the dynamics near different kinds of heteroclinic cycles and networks.


\subsection{Chirality}
\label{chirality_def}
Given a Bykov cycle, there are two different possibilities for the geometry of the flow around $\Gamma$,
depending on the direction trajectories turn around the one-dimensional heteroclinic connection from $O_1$ to $O_2$.
To make this rigorous, we need the following concepts adapted from \cite{LR2015}. 

\medbreak
Let $V_1$
and $V_2$ be small disjoint neighbourhoods of $O_1$ and $O_2$ with disjoint boundaries $\partial V_1$ and $\partial V_2$, respectively. Trajectories starting at $\partial V_1$ near $W^s(O_1)$ go into the interior of $V_1$ in positive time, then follow the connection from $O_1$ to $O_2$, go inside $V_2$, and then come out at $\partial V_2$. Let $\mathcal{Q}$ be a piece of trajectory like this from $\partial V_1$ to $\partial V_2$.
Now join its starting point to its end point by a line segment, forming a closed curve, that we call the  \emph{loop} of $\mathcal{Q}$.
The loop of $\mathcal{Q}$ and the cycle $\Gamma$ are disjoint closed sets.

We say that the two saddle-foci $O_1$ and $O_2$ in $\Gamma$ have the same \emph{chirality} if the loop of every trajectory is linked to $\Gamma$ in the sense that the two closed sets cannot be disconnected by an isotopy. Otherwise, we say that $O_1$ and $O_2$ have different chirality.

\subsection{Saturated horseshoe}
\label{ss: saturated horseshoe}
Given $(\tau_1, \tau_2) \in [0, \varepsilon]^2$, suppose that there is a cross-section $\mathcal{S}_\lambda$ to the flow $\varphi_{(\tau_1, \tau_2)}$ such that $\mathcal{S}_{(\tau_1, \tau_2)}$ contains a compact set $\mathcal{K}_{(\tau_1, \tau_2)}$ invariant by the first return map $\mathcal{F}_{(\tau_1, \tau_2)}$ to $\mathcal{S}_{(\tau_1, \tau_2)}$. Assume also that $\mathcal{F}_{(\tau_1, \tau_2)}$ restricted to $\mathcal{K}_{(\tau_1, \tau_2)}$ is conjugate to a full shift on a finite alphabet. Then the \emph{saturated horseshoe associated to $\mathcal{K}_{(\tau_1, \tau_2)}$} is the flow-invariant set $$\widetilde{\mathcal{K}_{(\tau_1, \tau_2)}}=\{\varphi_{(\tau_1, \tau_2)}(t,x)\,:\,t\in\RR,\,x\in \mathcal{K}_{(\tau_1, \tau_2)}\}.$$

\subsection{Rotational horseshoe}
\label{ss: rotational horseshoe}
Let $\mathcal{H} $ stand for the infinite annulus $\mathcal{H} = \EU^1 \times \RR$. We denote by $Homeo^+(\mathcal{H} )$ the set of homeomorphisms of the annulus which preserve orientation.
Given a homeomorphism $f :X \rightarrow X$  and a partition of $m>1$ elements $R_0,..., R_m$ of $X$, the itinerary  function $\xi: X \rightarrow \{0, ..., m-1\}^\ZZ= \Sigma_m$ is defined by $\xi(x)(j)=k$ if and only if $f^j(x)\in R_k$ for every $j\in \ZZ$.  
Following \cite{Passeggi}, we say that a compact invariant set $\Lambda \subset \mathcal{H} $ of $f \in Homeo^+(\mathcal{H} )$ is a rotational horseshoe if it admits a finite partition $P =\{R_0, ..., R_{m-1} \}$ with $R_i$ open sets of $\Lambda$ so that
\begin{enumerate}
\item the itinerary $\xi$ defines a semi-conjugacy between $f|_\Lambda$ and the full-shift $\sigma: \Sigma_m \rightarrow \Sigma_m$, that is $\xi  \circ f = \sigma \circ \xi$ with $\xi$ continuous and onto;
\item for any lift $F$ of $f$, there exist a positive constant $k$ and $m$ vectors $v_0, ...,v_{m-1} \in \ZZ \times \{0\}$ so that:
$$
\left\| (F^n(x)-x)  - \sum_{i=0}^n v_{\xi(x)}\right\| <k \qquad \text{for every} \qquad  x\in \pi^{-1}(\Lambda), \quad n\in \NN.
$$
\end{enumerate}

\subsection{SRB measure}
\label{app:SRB}
Given an attracting set ${\Omega}$ for a continuous map $R: \mathcal{M} \rightarrow \mathcal{M}$ of a compact manifold $ \mathcal{M}$, consider the Birkhoff average with respect to the continuous function $T:  \mathcal{M} \rightarrow \RR$ on the $R$-orbit starting at $x\in  \mathcal{M}$:
\begin{equation}
\label{limit1}
L(T, x)=\lim_{n\in \NN} \quad \frac{1}{n} \sum_{i=0}^{n-1} T \circ R^i(x).
\end{equation}

Suppose that, for Lebesgue almost all points $x\in \textbf{B}({\Omega})$,  the limit \eqref{limit1} exists and is independent on $x$. Then $L$ is a continuous linear functional in the set of continuos maps from  $\mathcal{M}$ to $\RR$ (denoted by $C(\mathcal{M}, \RR)$). By the Riesz Representation Theorem, it defines a unique probability measure $\mu$ such that:
\begin{equation}
\label{limit2}
\lim_{n\in \NN} \quad \frac{1}{n} \sum_{i=0}^{n-1} T \circ R^i(x) = \int_{\Omega} T \, d\mu
\end{equation}
for all $T\in C(\mathcal{M}, \RR)$ and for Lebesgue almost all points $x\in \textbf{B}({\Omega})$.  If there exists an ergodic measure $\mu$ supported in ${\Omega}$ such that \eqref{limit2} is satisfied for all continuous maps $T\in C(\mathcal{M}, \RR)$ for Lebesgue almost all points $x\in \textbf{B}({\Omega})$, where $\textbf{B}({\Omega})$ has positive Lebesgue measure, then $\mu$ is called a SRB measure and ${\Omega}$ is a SRB attractor.

\subsection{Symmetry and lifting by rotation}
\label{lift1}

Given a group $\mathcal{G}$ of endomorphisms of $\EU^3$, we will consider two-parameter families of vector fields $(f_{(\tau_1, \tau_2)})$ under the equivariance assumption $$f_{(\tau_1, \tau_2 )}(\gamma x)=\gamma f_{(\tau_1, \tau_2 )}(x)$$ for all $x \in \EU^3$, $\gamma \in \mathcal{G}$ and $(\tau_1, \tau_2 )\in  [0, \varepsilon]^2.$
For an isotropy subgroup $\widetilde{\mathcal{G}}< \mathcal{G}$, we will write $\Fix(\widetilde{\mathcal{G}})$ for the vector subspace of points that are fixed by the elements of $\widetilde{\mathcal{G}}$. Observe that, for $\mathcal{G}-$equivariant differential equations, the subspace $\Fix(\widetilde{\mathcal{G}})$ is flow-invariant.

The authors of \cite{ACL06, Rodrigues} investigate how some properties of a $\ZZ_2$--equivariant vector field on $\RR^n$ lift by a rotation to properties of a corresponding vector field on $\RR^{n+1}$. For the sake of completeness,  we review some of these properties.

\bigbreak
 Let $f_n$ be a $\ZZ_2(\gamma_n)$--equivariant vector field on $\RR^n$. Without loss of generality, we may assume that $f_n$ is equivariant by the action of
$$\gamma_n(x_1, x_2, ...., x_{n-1}, y)=  (x_1, x_2, ...., x_{n-1}, - y).$$
The vector field $f_{n+1}$ on $\RR^{n+1}$ is obtained by adding the auxiliary equation $\dot\theta=\omega>0$ and interpreting $(y, \theta)$ as polar coordinates. In cartesian coordinates $(x_1, . . . , x_{n-1}, r_1, r_2) \in \RR^{n+1}$, this equation corresponds to the system $r_1 = |y| \cos \theta$ and $r_2 = |y| \sin \theta$. The resulting vector field $f_{n+1}$ on $\RR^{n+1}$ is called the \emph{lift by rotation of} $f_n$, and is $\mathbb{SO}(2)$--equivariant in the last two coordinates.
\medbreak

Given a set $\Lambda \subset \RR^n$, let $\mathcal{L}(\Lambda)\subset \RR^{n+1}$ be the lift by rotation of $\Lambda$, that is,
$$\Big\{(x_1,...,x_{n-1}, r_1, r_2) \in \RR^{n+1} \colon \, (x_1,\ldots,x_{n-1},||(r_1,r_2)||) \quad \text{or} \quad (x_1,\ldots,x_{n-1},-||(r_1,r_2)||) \in \Lambda\Big\}.$$
It was shown in \cite[Section 3]{ACL06} that, if $f_n$ is a $\ZZ_2(\gamma_n)$--equivariant vector field in $\RR^n$ and $f_{n+1}$ is its lift by rotation to $\RR^{n+1}$, then:
\medskip
\begin{enumerate}
\item If $p$ is a hyperbolic equilibrium of $f_n$ outside $\text{Fix}(\ZZ_2(\gamma_n))$, then $\mathcal{L}(\{p\})$ is a hyperbolic periodic solution of $f_{n+1}$ with minimal period $2\pi/\omega$.
\medskip
\item If $p$ is a hyperbolic equilibrium of $f_n$ lying in $\text{Fix}(\ZZ_2(\gamma_n))$, then $\mathcal{L}(\{p\})$ is a hyperbolic equilibrium of $f_{n+1}$.
\medskip
\item If $[p_1 \to p_2]$ is a $k$-dimensional heteroclinic connection between equilibria $p_1$ and $p_2$ and it is not contained in $\text{Fix}(\ZZ_2(\gamma_n))$, then it lifts to a $(k + 1)$--dimensional connection between the periodic orbits $\mathcal{L}(\{p_1\})$ and $\mathcal{L}(\{p_2\})$ of $f_{n+1}$.
\medskip
\item If $\Lambda$ is a compact $f_n$--invariant asymptotically stable set, then $\mathcal{L}(\Lambda)$ is a compact $f_{n+1}$--invariant
asymptotically stable set.

\end{enumerate}

\end{document}